\documentclass{paper}

\usepackage{arxiv}

\usepackage[utf8]{inputenc} 
\usepackage[T1]{fontenc}    
\usepackage{hyperref}       
\usepackage{url}            
\usepackage{booktabs}       
\usepackage{amsfonts}       
\usepackage{nicefrac}       
\usepackage{microtype}      
\usepackage{subcaption}
\usepackage{mathrsfs,bm,siunitx}
\usepackage{mathtools}
\usepackage{ragged2e}
\usepackage[]{graphicx}
\usepackage[]{color}
\usepackage{alltt}
\usepackage{mathtools}
\usepackage{amssymb}
\usepackage{aligned-overset}
\usepackage{amsmath}
\numberwithin{equation}{section}
\usepackage{tikz}
\usepackage{capt-of}
\usepackage{xcolor,pgf}
\usepackage{color}
\usepackage{transparent,mathtools}
\usepackage{caption}
\usetikzlibrary{arrows}
\usetikzlibrary{decorations.markings}
\usepackage{pgfplots}
    \usetikzlibrary{
        angles,
        quotes,
    }
\usepackage{comment}
\usepackage[square, comma, numbers, sort]{natbib}

\usepackage{lipsum}
\usepackage{graphicx}
\usepackage{amsmath,amsthm}
\numberwithin{equation}{section}
\usepackage{color}
\usepackage{colortbl}
\usepackage{xcolor}
\colorlet{LighterGray}{white!70!gray}
\colorlet{LightererGray}{white!85!gray}
\usepackage{transparent}
\usepackage{tikz}
\usetikzlibrary{arrows}
\usetikzlibrary{decorations.markings}
\usepackage{pgfplots}
\pgfplotsset{compat=1.14}
    \usetikzlibrary{
        angles,
        quotes,
    }
    
    \usepackage{multirow}
\usepackage[utf8]{inputenc}
\usepackage{ragged2e}
\usepackage{xcolor}
\colorlet{LighterGray}{white!70!gray}
\colorlet{LightererGray}{white!85!gray}
\usepackage[square, comma, numbers, sort]{natbib}
\usepackage{fancyhdr}
\newcommand{\changefont}{%
    \fontsize{8}{11}\selectfont
}
\usepackage{hyperref}
\hypersetup{
  colorlinks = true,
  linkcolor = magenta,
  urlcolor = magenta,
  citecolor = cyan}
\usepackage{scalerel,stackengine,amsmath}

\usepackage{enumitem}
\newlist{abbrv}{itemize}{1}
\setlist[abbrv,1]{label=,labelwidth=1in,align=parleft,itemsep=0.1\baselineskip,leftmargin=!}
\usepackage{listings}
\usepackage{color} 
\definecolor{mygreen}{RGB}{28,172,0} 
\definecolor{mylilas}{RGB}{170,55,241}

\usepackage{amstext} 
\usepackage{array}   
\newcolumntype{L}{>{$}l<{$}}
\usepackage{bm}
\DeclareMathOperator*{\argmin}{arg\,min}
\usepackage{pgfplots}
\pgfplotsset{compat=1.14}
\theoremstyle{definition}

\title{Development of a photothermal measurement model to determine layer thickness of multi-layered coating systems with unknown thermal properties}

\author{
  Dimitri Rothermel \\
  Department of Mathematics\\
  Saarland University\\
  Saarbr\"ucken, Germany \\
  \texttt{rothermel@math.uni-sb.de} \\
       \And
  Thomas Schuster \\
  Department of Mathematics\\
  Saarland University\\
  Saarbr\"ucken, Germany \\
  \texttt{thomas.schuster@num.uni-sb.de} \\
}

\begin{document}
\maketitle

\begin{abstract}
In this article, a general model for 1D thermal wave interference is derived for multi-layered coating systems on a thermally thick substrate using the same principles as for the well established one-layered and two-layered coating cases. Using the lock-in thermography principle, an illumination source modulates the surface of those systems periodically by a planar, sinusoidal wave form with a fixed frequency. The coating systems absorb the optical energy on its surface and convert it into thermal energy, resulting in the propagation of a spatially and temporally periodic thermal wave with the same frequency. These thermal waves, originating at the surface, are reflected and transmitted at each interface leading to infinitely many wave trains that need to be tracked in order to formulate the final surface temperature as a superposition of all these waves. The heat transfer inside the object depends on the layer thickness of each coating, but also on the thermal properties of each layer material. The goal is to have a mathematical and physical model which describes the phase angle data measured by an infrared camera. Having these data, the main objective of this paper is to determine the thickness of each coating layer. In practice, the thermal properties of the layers usually are unknown, which makes the process even more difficult. For that reason, this article presents a concept to determine the thermal properties in advance.

\end{abstract}



\section{Introduction}\label{sec:intro}


In this aticle we investigate general multi-layered coating systems of $n+1$ total layers, where $n\in \mathbb{N}$ coating layers of different materials $M_1,\dots,M_n$ are applied on top of each other on a substrate material $M_{n+1},$ see Figure \ref{fig:multi-layered_systems}.

\begin{center}
\captionsetup{type=figure}
 \resizebox{8cm}{!}{\begin{tikzpicture}[
    scale=5,
    axis/.style={very thick, >=stealth'},
    important line/.style={thick},
    dashed line/.style={dashed, thin},
    pile/.style={thick, ->, >=stealth', shorten <=2pt, shorten
    >=2pt},
    every node/.style={color=black}
    ]
    \draw[white] (-1,1.2)--(-1,1.4)--(1,1.4)--(1,1.2)--(-1,1.2) --cycle;
    \draw[fill=lightgray] (-1,0)--(-1,0.6)--(1,0.6)--(1,0)--(-1,0) --cycle;
     \draw[fill=LighterGray] (-1,0.6)--(-1,1.2)--(1,1.2)--(1,0.6)--(-1,0.6) --cycle;
      \draw[fill=LightererGray] (-1,1.4)--(-1,2)--(1,2)--(1,1.4)--(-1,1.4) --cycle;
    \draw[fill=gray] (-1,0)--(1,0)--(1,-0.75)--(-1,-0.75)--(-1,0) --cycle;
    

    \draw[very thick,<->,>=stealth] (1.1,0) -- (1.1,0.6) node[xshift=1.0cm,yshift=-1.5cm] {\huge{$\bm{L_{n}}$}};
    \draw[very thick,<->,>=stealth] (1.1,0) -- (1.1,-0.75) node[xshift=1.2cm,yshift=2cm] {\huge{$\bm{L_{n+1}}$}};
    \draw[very thick,<->,>=stealth] (1.1,0.6) -- (1.1,1.2) node[xshift=1.3cm,yshift=-1.5cm] {\huge{$\bm{L_{n-1}}$}};
    \draw[very thick,<->,>=stealth] (1.1,1.4) -- (1.1,2) node[xshift=1.0cm,yshift=-1.5cm] {\huge{$\bm{L_{1}}$}};

    \draw[axis] (-1,0)  -- (1,0) node[above, xshift=-3cm, yshift = 1cm] {\huge{$\bm{M_{n}(\bm{\alpha}_{n},\bm{e}_{n})}$}} node[below, xshift=-4.1cm,yshift=-1cm] {\huge{$\bm{M_{n+1}(\bm{\alpha}_{n+1},\bm{e}_{n+1})}$}} node[above, xshift=-4.1cm,yshift=4cm] {\huge{$\bm{M_{n-1}(\bm{\alpha}_{n-1},\bm{e}_{n-1})}$}} node[above, xshift=-3cm,yshift=8cm] {\huge{$\bm{M_1(\bm{\alpha}_1,\bm{e}_1)}$}} node[above, xshift=-5cm,yshift=6.25cm] {$\bm{\vdots}$};
  \end{tikzpicture}}
  
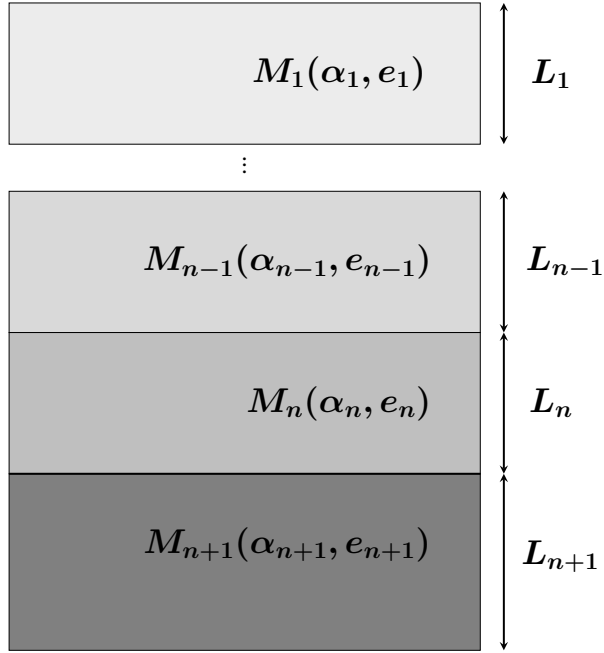
\captionof{figure}{System of $n$ coating layer materials $M_1$,\dots, $M_{n}$ on a substrate material $M_{n+1}$}
 \label{fig:multi-layered_systems}
  \end{center}

For $i=1,\dots,n+1$, the variable $L_i \in \mathbb{R_+}$ denotes the thickness of the associated layer. The thermal properties of each material are characterized by the thermal \emph{diffusivity} and the \emph{thermal effusivity} and denoted by $\alpha_i\in \mathbb{R_+}$ and $e_i\in \mathbb{R_+}$, respectively.

The determination of the vector
$$ \mathbf{L}:=(L_1,\dots,L_n)^T \in \mathbb{R}_+^n$$
having the coating layer thicknesses as components
by a nondestructive testing method is of great interest in the manufacturing process and quality control of such systems. In this paper, it is suggested that photothermal methods, such as infrared thermography, are specially suited for this purpose. A planar and periodically modulated light source is used to irradiate the top surface of a test object, which leads to a propagating thermal wave inside the object. Such thermal waves show the same behavior at interfaces as shear waves in elasticity theory or acoustic and optical waves in the visible spectral range, i.e. reflection and transmission coefficients can be calculated. The temperature response, which can be measured by an infrared camera, is a result of thermal wave interference, i.e. the superposition of thermal wave trains, which of course depend on the coating layer thicknesses and its material thermal properties. In our setup the phase angle of the temperature response is measured.

Unfortunately, the thermal properties
$$(\alpha_i,e_i) \text{, \ \ \ } i=1,\dots,n+1$$
are often partially or even completely unknown making the process of determining the coating layer vector $\mathbf{L}\in \mathbb{R}_+^n$ even more difficult. The goal of this paper is to develop a concept utilizing the generated data of multiple samples in a way, where the unknown thermal properties can be determined as well in order to calibrate a model for the layer thickness determination.

\section{Motivation}

The determination of thicknesses even in one-layered coating systems is of great interest for quality control, e.g. in the manufacturing process of electrode coatings of Lithium-ion batteries, where a Lithium cobalt oxide coating of $50$ to $100\ \mu m$ is applied to a thin ($\pm 10\ \mu m$) aluminium substrate. Undesired coating thicknesses lead to performance issues of the battery or higher production costs. In the worst case, dangerous situations can arise, such as the explosion of batteries, cf. \cite{wang2012thermal}. Another application example where the control of layer thicknesses is very important, is the application of non-electrolytic zinc-aluminium flakes to protect metallic surfaces by increasing its corrosion resistance, cf. \cite{fourez1993application}, \cite{hare1982anti}. A further application example showing a higher complexity, i.e. with two coating layers, is the plastic housing of laptops or smartphones. For example, the substrate might be some rigid ABS plastic of $1$ to $2\ mm$, the first thin coating layer (the so-called basecoat) might be some polyurethane ($\pm 20 \ \mu m$) that acts like a thermal insulator against the rising heat of the electronics, and the second coating layer (the so-called topcoat of $\pm 20\ \mu m$) might be some UV hardened resin acting like a visually appealing surface finish. To give a last idea of an application with increasing coating layers, it is certainly worth mentioning the automotive and aviation industry (\cite{ellrich2020terahertz}, \cite{krimi2016highly}), where metallic or CFRP substrates are coated with different paints, e.g. a primer, basecoat, topcoat and clearcoat. Here again, each layer fulfills its own function and it is therefore necessary to control the individual layer thicknesses. Of course, like for every other industry sector with an economical objective, the goal is also to reduce resources, while maintaining the full functionality of the product. 

Obviously, depending on the manufacturing process, discrepancies between target and actual thickness values can arise. It is therefore necessary to build devices and develop algorithms in order to keep track of those unwanted deviations as early as possible. Instead of the conservative approach of performing quality control randomly after the completion of a bigger production line by cutting up individual samples to evaluate cross-sectional slice images under the microscope, the modern Industry 4.0 (cf. \cite{lasi2014industry}) approach requires inline procedures that are nondestructive\footnote{In general, nondestructive testing or evaluation summarizes numerous techniques that aim to test and evaluate the properties of materials without causing damage.} and monitor the process in real time. 

There are a few things to consider when choosing a measuring process (cf. \cite{gholizadeh2016review}) in this special setting. Eddy current or inductive based measurement methods only work on metallic substrates, X-ray or beta backscattering methods only work with certain metal groups and also require compliance with strict occupational health and safety, radiation protection and disposal measures. Furthermore, ultrasonic and capacitive methods need contact with the test specimen and are therefore not suitable for measuring wet coatings (e.g. for the automotive paint process line) or uncured powder coatings. A very promising method, which is contact-free and uses harmless and non-invasive electromagnetic radiation, is presented by the terahertz technology, see \cite{krimi2016highly}. Unfortunately, the usage of highly sensitive devices is needed, such as a femtosecond laser (which is able to keep track of signals in the picosecond range and) which is often prohibitively expensive for medium-sized companies.

In this work, it is suggested that \textit{photothermal methods}, especially \textit{infrared thermography}, could solve the aforementioned problems. Such methods are contact-free with the test object and have manageable costs. A typical setting is presented in Figure \ref{fig:setting},

\begin{center}
\captionsetup{type=figure}
    \includegraphics[width = 0.9\textwidth]{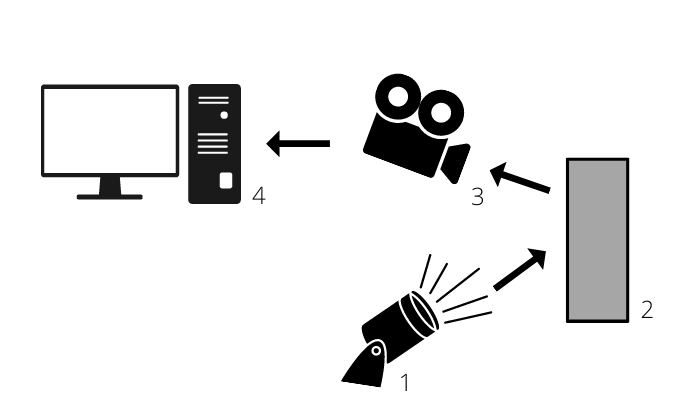}
\captionof{figure}{Measurement principle of infrared thermography} \label{fig:setting}
\end{center}

where the components and the operating principle are as follows:

\begin{itemize}
    \item[1] \textit{Optical excitation sources}, such as lamps or lasers (cf. \cite[Chapter 4]{ALM1996}), irradiate an object under investigation with optical (i.e. visible) light. Its task is to produce the \textit{photothermal effect}, i.e.
    \item[2] the \textit{test object} absorbs the optical energy (light) and converts it into thermal energy (heat). Heat transfer inside the test object occurs, which is dependent on the characteristics of the object geometry, its material composition and thermal properties. 
    \item[3] In the next step, the thermal response of the test object is recorded by a \textit{measuring device}, such as an infrared camera.
    \item[4] Finally, a \textit{computer} processes and evaluates the data to provide certain properties (or even defects) of the considered test object. 
\end{itemize}

In particular, it should be mentioned that the photothermal method requires only two conditions:
\begin{itemize}
    \item The investigated coating must be susceptible to optical radiation of a certain wavelength range, i.e. for near infrared, visible or UV light.
     \item There must be a thermal contrast between two adjacent layers. Otherwise, a transparent interface would be created without significant reflection.
\end{itemize}

As a brief note, there are mainly two classical optical excitation types in thermographic processes, i.e. the lock-in thermography and pulsed thermography, see e.g. \cite{wu1998lock}, \cite{breitenstein2003lock}, \cite{chatterjee2011comparison}.  The latter method analyses the transient response and propagation of heat pulses in the test object. A prominent processing technique in pulsed thermography is the TSR (Thermographic Signal Reconstruction) method, where the logarithmic time derivatives of the signal are examined, cf. \cite{shepard2003reconstruction}, \cite{shepard2015advances}, \cite{balageas2015thermographic}. In lock-in thermography, the surface of the test object is periodically modulated by a planar sinusoidal wave form with a fixed frequency. Radiation is absorbed and leads to the propagation of a spatially and temporally periodic temperature field in the test object, which is also referred to as the so-called \textit{thermal wave}. Since this thermal wave has the same frequency as the excitation, especially the phase angle (or phase difference) carries a lot of information. Thermal waves show the same behavior at interfaces as shear waves in elasticity theory or acoustic and optical waves in the visible spectral range, i.e. reflection and diffraction occur. Therefore, thermal reflection and transmission coefficients can be derived, see e.g. \cite{greiner2008klassische}, \cite{macke1962thermodynamik}. 

In this paper, the  frequency lock-in principle is used and the test object is a multi-layered coating system, see Figure \ref{fig:multi-layered_systems}. The parameter estimation of the vector $\mathbf{L}\in \mathbb{R}_+^n$ is very suitable to be understood in terms of inverse problems (cf. \cite{louis2013inverse}, \cite{schuster}, \cite{rothermel2021determination}, \cite{rothermel2021solving}), since interior properties that are not directly accessible are to be determined by processing exterior data in the form of the phase angle of the measured surface temperature.

\section{Mathematical Setting}

In order to understand the physical process of thermal wave interference for the general case of multi-layered coating systems, it is necessary to take a look at known models for the semi-infinite medium and for the cases with $n=1$ and $n=2$ coatings. The following mathematical notations and formulations in the next subsections can be found in \cite{CAR1947},\cite{palmer2010art}, \cite{ALM1996}, \cite{lepoutre1985nondestructive}.

\subsection{Basics of thermal wave generation}
\label{wavegeneration}

Thermal waves can be mathematically characterized as solutions of the heat diffusion equation. The type of the heat source at the surface, which represents the appropriate boundary conditions, influences the surface temperature distribution and determines the generation of waves. The most common type of excitation is a periodic, planar energy input by high-performance laser beams of a single specific excitation frequency, i.e., lock-in excitation. 
\newline
\newline
For the sake of simplicity, we firstly consider an isotropic homogeneous semi-infinite medium $M$ (that means an infinite extension of the medium in $x$-direction), which is described by the following figure: 
\begin{center}
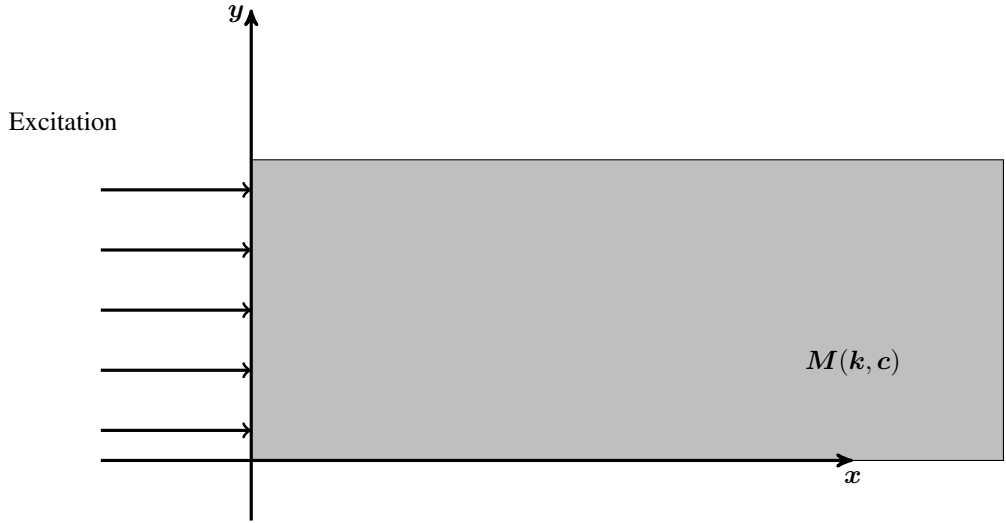

\captionsetup{type=figure}
 \begin{tikzpicture}[
     scale=4.,
    axis/.style={very thick, >=stealth'},
    important line/.style={thick},
    dashed line/.style={dashed, thin},
    pile/.style={thick, ->, >=stealth', shorten <=2pt, shorten
    >=2pt},
    every node/.style={color=black}
    ]
     \draw[fill=lightgray] (0,0)--(0,1)--(2.5,1)--(2.5,0)--(0,0) --cycle node[above, xshift=8cm,yshift=1cm]{$\bm{M}(\bm{k},\bm{c})$};
    \draw[very thick, ->] (-0.5,0.5) -- (0,0.5); 
    \draw[very thick, ->] (-0.5,0.7) -- (0,0.7);
    \draw[very thick, ->] (-0.5,0.3) -- (0,0.3);
    \draw[very thick, ->] (-0.5,0.1) -- (0,0.1);
    \draw[very thick, ->] (-0.5,0.9) -- (0,0.9);
       \draw[axis,->] (-0.5,0)  -- (2,0)  node[below]{$\bm{x}$};
    \draw[axis,->] (0,-0.2) -- (0,1.5) node[above,xshift=-0.2cm,yshift=-0.3cm]{$\bm{y}$}
     node[above, xshift=-2.5cm,yshift=-1.75cm]{Excitation};
\end{tikzpicture}
\captionof{figure}{Thermal wave generation and propagation in a semi-infinite medium}
\end{center}

Here, $k$ denotes the thermal \emph{conductivity} and $c$ denotes the volumetric \emph{heat capacity} of the material $M$. The thermal diffusivity can then be calculated by $\alpha:=\frac{k}{c}$ and the thermal effusivity by $e:=\sqrt{kc}$.

We assume that the heated surface occupies the $y-z$-plane at $x=0$. Consequently, to obtain the temperature distribution at the surface of the medium, we must solve the Fourier equation, which in this situation reduces to the one-dimensional case:
\begin{align}\label{eq1}
    \frac{\partial ^2 T}{\partial x^2} = \frac{1}{\alpha} \frac{\partial T}{\partial t}, \quad x,t>0.
\end{align}
We at first need to specify the boundary conditions. 
We excite the medium's surface by a plane harmonic, thus temporal heating of modulation frequency $\omega\coloneqq 2\pi f$ for some frequency $f$ and source intensity $Q_0$. This can be described by an excitation term of the form
\begin{align*}
    \frac{Q_0}{2}\left[1+\cos(\omega t)\right],
\end{align*}
yielding the generation of thermal waves in the inside of the medium.
\newline
\newline
Since the periodic thermal energy is subject to conduction into the solid, by using the appropriate rate equation, we obtain the following boundary condition on the surface of the medium: 
\begin{eqnarray}
    -\kappa \frac{\partial T}{\partial x}&=&\text{Re} \left\{\frac{Q_0}{2}\left[1+\exp{(i\omega t)}\right]\right\}\nonumber\\
    &=&\frac{Q_0}{2}\left[1+\cos(\omega t)\right]\nonumber\\
    \label{BC}
    &=&\underbrace{\frac{Q_0}{2}}_{\text{dc component}} + \ \ \ \ \  \underbrace{\frac{Q_0}{2}\cos(\omega t)}_{\text{ac component}} , \quad x=0,t>0.
\end{eqnarray}
Here, DC means the \emph{Direct Current} and AC the \emph{Alternating Current}.
Neglecting the dc component for a moment, since this quantity will not be relevant in the later applications, and applying a time-harmonic approach yield
\begin{align*}
    T(x,t) = \text{Re}\left[T(x)\exp(i\omega t)\right].
\end{align*}
Plugging this into Equation \eqref{eq1}, we end up with 
\begin{align*}
    \exp(i\omega t)\left( \frac{\partial ^2 T(x)}{\partial x^2}-\frac{i \omega}{\alpha}T(x)\right) = 0.
\end{align*}
Taking into account that $T(x)$ must be finite for $x\to+\infty$, we receive the solution of the boundary value problem \eqref{eq1}, \eqref{BC} as
\begin{align}
    T(x,t) = \frac{Q_0}{2\kappa \sigma} \exp(-\sigma x + i\omega t), \quad \sigma\coloneqq (1+i) \sqrt{\frac{\omega}{2\alpha}}.
\end{align}
By a multiplication with $1=\frac{i+1}{\sqrt{2}}\exp{(-i\frac{\pi}{4})}$ and further simplifications, we obtain a more significant expression given as
\begin{align}
   \bm{T}(x,t) = \frac{Q_0}{2e\sqrt{\omega}}\exp\left(-\frac{x}{\mu}\right)\exp{\left[i\left(\omega t -\frac{x}{\mu}-\frac{\pi}{4}\right)\right]},
\end{align}
where \begin{align}
    \mu \coloneqq \sqrt{\frac{2\alpha}{\omega}}
\end{align}
is the so-called \textit{thermal diffusion length}.

Hence, thermal waves are significantly damped and $\mu$ controls the penetration depth into the material. For small thermal diffusivity $\alpha$ the thermal waves do only slightly propagate into the interior of the material. In contrast, by decreasing the modulation frequency $\omega$, we obtain a deeper penetration of the thermal waves into the material. This phenomenon is very useful in the photothermal measurement of layer thicknesses.

Furthermore, note that there occurs a progressive phase shift by
\begin{align}\label{phaseshift}
    \varphi = -\frac{x}{\mu} - \frac{\pi}{4}. 
\end{align}
between the temperature at the surface and a point $x$ located at
the propagating thermal wave in the material. Thus, at the surface $x=0$ there is an expressive phase difference of $-45$ degree between the excitation source and the resulting surface temperature. 

\subsection{Transmission and reflection}

When the irradiated object is not a semi-infinite medium but a composition of at least two materials $M$ and $\Tilde{M}$ (with thermal effusivities $e$ and $\Tilde{e}$, respectively), the thermal wave travels trough $M$ first towards $\Tilde{M}$ and when the planar thermal wave propagation direction is perpendicular to its interface, the thermal reflection and transmission coefficients are
\begin{align}
    R = \frac{1-b}{1+b}, \quad T=1+R=\frac{2}{1+b},
\end{align}

where 
\begin{align}
    \bm{b} = \frac{\Tilde{e}}{e}
\end{align}
is the so-called \textit{thermal refraction index}, which characterizes the thermal contrast between the two media. If there is no thermal contrast, i.e. when $e=\Tilde{e}$, it follows that $R=0$. In that case there would be no significant reflection from that interface and therefore no contribution to thermal wave interference effects influencing the surface temperature. Thus, regarding the determination of coating layer thicknesses it is crucial to guarantee that materials are distinguishable.

In the following, whenever we add indices $j=1,\dots,n$ to the reflection or transmission coefficient, we refer to the interface between the materials $M_j$ and $M_{j+1}.$ The direction of these coefficients has to be understood downwards, cf. Figure \ref{fig:multi-layered_systems}. For the upward direction, i.e. when the thermal wave approaches from below, we add a prime notation.  For example,
\begin{align*}
    R_n &:= R_{M_{n}\to M_{n+1}},\\
    T_1' &:= T_{M_2 \to M_1}.
\end{align*}

The only exception regarding the direction is for $R_0$ and $T_0$, where the top surface of first layer material $M_1$ is exposed to air $M_0$,

i.e.
\begin{align*}
    R_0 &:= R_{M_{1}\to M_{0}},\\
    T_0 &:= T_{M_1 \to M_0}.
\end{align*}

We refer the reader to \cite{ALM1996} for a more detailed derivation of the above expressions.

\subsection{Basics of thermal wave interference}

In this subsection, we want to discuss the basics of thermal wave interference by investigating the cases of $n=1$ and $n=2$ coatings. Here, mathematical formulas are well established, cf. \cite{ALM1996}, \cite{palmer2010art} and \cite{lepoutre1985nondestructive}. For the convenience of the reader, we summarize the findings before we extend the insights to multi-layered coating systems in the next subsection.

\subsubsection{One-layered coatings on a substrate}

Consider the following system of two layers consisting of media $M_1$ and $M_2$, possessing different thermophysical properties: 
\begin{center}
\captionsetup{type=figure}
 \begin{tikzpicture}[
    scale=5,
    axis/.style={very thick, >=stealth'},
    important line/.style={thick},
    dashed line/.style={dashed, thin},
    pile/.style={thick, ->, >=stealth', shorten <=2pt, shorten
    >=2pt},
    every node/.style={color=black}
    ]
    \draw[fill=lightgray] (-1,0)--(-1,0.6)--(1,0.6)--(1,0)--(-1,0) --cycle;
    
    \draw[fill=gray] (-1,0)--(1,0)--(1,-0.75)--(-1,-0.75)--(-1,0) --cycle;
    
    \draw[->,>=stealth] (0,0.6) -- (0,1)node[above,yshift=1cm] {Transmitted thermal waves};
    \draw[->,>=stealth] (0.5,0.6) -- (0.5,1) ;
    \draw[ ->,>=stealth] (-0.5,0.6) -- (-0.5,1);
    \draw[->,>=stealth] (0,0) -- (0,-0.5)node[above,yshift=-1cm] {Transmitted thermal waves} ;
    \draw[->,>=stealth] (0.5,0) -- (0.5,-0.5) ;
    \draw[->,>=stealth] (-0.5,0) -- (-0.5,-0.5);
    \draw[->,>=stealth] (-0.5,0.55) -- (-0.52,0) -- (-0.54,0.6) -- (-0.56,0) -- (-0.58,0.6) -- (-0.6,0) -- (-0.62,0.6);
    \draw[->,>=stealth] (-0.1,0.55)  -- (-0.12,0) -- (-0.14,0.6) -- (-0.16,0) -- (-0.18,0.6);
    \draw[->,>=stealth] (0.1,0.55)  -- (0.12,0.6) -- (0.14,0) -- (0.16,0.6) --(0.18,0) -- (0.2,0.6);
    \draw[->,>=stealth] (0.5,0.55) -- (0.52,0.6)  -- (0.54,0) -- (0.56,0.6) -- (0.58,0) -- (0.6,0.6);
    \draw[->,>=stealth] (-0.45,0.55) -- (-0.45,0.6);
    \draw[->,>=stealth] (0.45,0.55) -- (0.45,0.6);
    \draw[very thick,<->,>=stealth] (1.1,0) -- (1.1,0.6) node[xshift=0.5cm,yshift=-1.5cm] {$\bm{L_1}$};
    \draw[very thick,<->,>=stealth] (1.1,0) -- (1.1,-0.75) node[xshift=0.5cm,yshift=2cm] {$\bm{L_2}$};
    \draw[axis] (-1,0)  -- (1,0) node[above, xshift=-1.5cm, yshift = 1cm] {$\bm{M_1(\alpha_1,e_1)}$} node[below, xshift=-1.5cm,yshift=-1cm] {$\bm{M_2(\alpha_2,c_2)}$} node[above, xshift=-1.5cm,yshift=4cm] {$\bm{M_0(\alpha_0,e_0)}$};
  \end{tikzpicture}
  
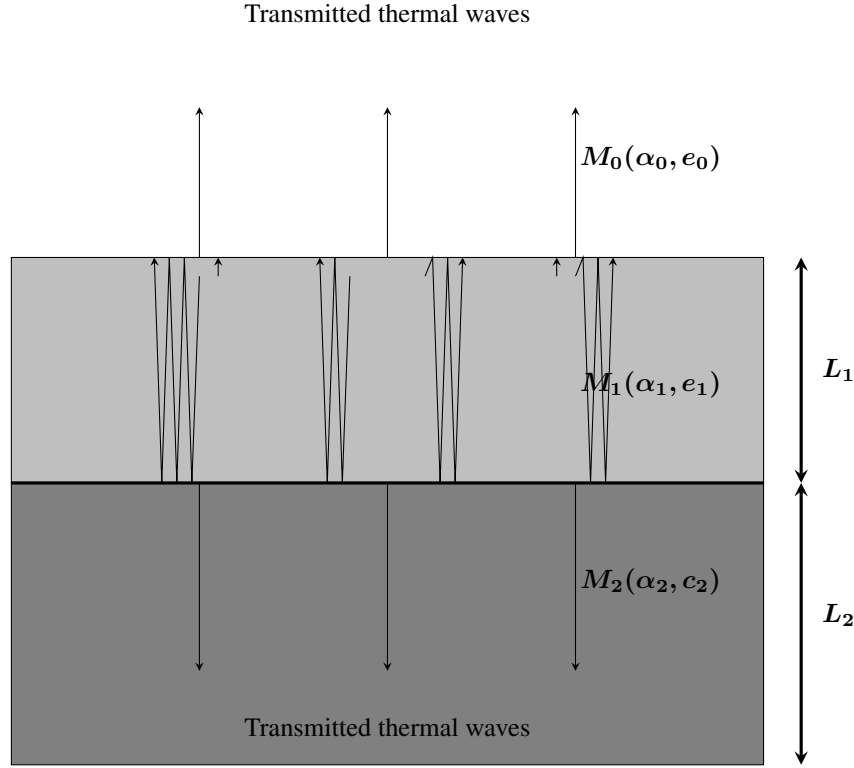
\captionof{figure}{Thermal wave interference in one coating layer on a thermally thick substrate material}
  \end{center}

Assume that both media $M_1$ and $M_2$ have homogeneous thermophysical properties and that $M_1$ has thickness $L_1$ and $M_2$ has thickness $L_2$. Furthermore, assume that $M_2$ is \textit{thermally thick}. That means that $L_2$ is significantly larger than the thermal diffusion length $\mu_2$. This assumption is not unusual as the substrate is often much thicker then the applied coating. This guarantees that the transmitted thermal waves into the substrate do not have an effect on the surface temperature because they can be neglected due to the large attenuation. Moreover, assume that the whole system is exposed to air, that we denote by $M_0$.
Suppose that the surface is illuminated by a plane, normal and periodic heating. 
\newline 
The single wave trains, that are generated near the surface $x=0$, propagate towards the interface between the two media and back towards the surface of $M_1$. When meeting any interfaces, the waves are partially reflected and transmitted.
Thermal wave interference effects occur, meaning that the surface temperature at $x=0$ is a sum of all thermal wave trains. In general, when a thermal wave has traveled a distance $x>0$ the amplitude will be damped by $\exp(-\sigma x)$, where $\sigma=(1+i)\frac{1}{\mu}$ is the complex wave number. Hence, for the first reflection order wave train in material $M_1$ the thermal wave has a propagation path of length $2L_1$ leading to an attenuation by $\exp(-2\sigma_1 L_1)$. We emphasize that the following figures include sketches of thermal wave trains and give the impression that there is a lateral diffusion in the material itself. This is not the case, since we have one-dimensional propagation. The sketches serve only for graphical visualization.

We note that we do not consider any bulk absorption of the absorbed radiation in this paper. This means that the photothermal effect provides a surface heating only, i.e. thermal waves originate exclusively near $x=0$. The presented literature sources distinguish between the following two cases:

\begin{itemize}
    \item[\textbf{(1)}] Waves that are first reflected from the interface between $M_0$ and $M_1$ as described by the following figure:
  \begin{center}
  \captionsetup{type=figure}
 \begin{tikzpicture}[
    scale=5,
    axis/.style={very thick, >=stealth'},
    important line/.style={thick},
    dashed line/.style={dashed, thin},
    pile/.style={thick, ->, >=stealth', shorten <=2pt, shorten
    >=2pt},
    every node/.style={color=black}
    ]
    \draw[fill=lightgray] (-1,0)--(-1,0.6)--(1,0.6)--(1,0)--(-1,0) --cycle;
    
    \draw[fill=gray] (-1,0)--(1,0)--(1,-0.75)--(-1,-0.75)--(-1,0) --cycle;
    \draw[->,>=stealth] (0.1,0.55)  -- (0.12,0.6) -- (0.14,0) -- (0.16,0.6) --(0.18,0) -- (0.2,0.6)node[above,xshift=-0.75cm,yshift=-1.5cm] {$\bm{a_2}$};
    \draw[->,>=stealth] (0.5,0.55) -- (0.52,0.6)  -- (0.54,0) -- (0.56,0.6) -- (0.58,0) -- (0.6,0.6) -- (0.62,0) -- (0.64,0.6) node[above,xshift=-1cm,yshift=-1.5cm] {$\bm{a_3}$} node[above,xshift=1.2cm,yshift=-1.5cm] {$\dots$};
    \draw[->,>=stealth] (-0.1,0.55)  -- (-0.12,0.6) -- (-0.14,0) -- (-0.16,0.6) --(-0.18,0) -- (-0.2,0.6) -- (-0.22,0) -- (-0.24,0.6) -- (-0.26,0) -- (-0.28,0.6)node[above,xshift=-0.5cm,yshift=-1.5cm] {$\bm{a_4}$};
    \draw[->,>=stealth] (-0.5,0.55) -- (-0.52,0.6)  -- (-0.54,0) -- (-0.56,0.6) -- (-0.58,0) -- (-0.6,0.6) -- (-0.62,0) -- (-0.64,0.6) -- (-0.66,0) -- (-0.68,0.6)-- (-0.7,0) -- (-0.72,0.6)--(-0.74,0) -- (-0.76,0.6)node[above,xshift=-0.5cm,yshift=-1.5cm] {$\bm{a_6}$};
     \draw[fill=gray] (-1,0)--(1,0)--(1,-0.75)--(-1,-0.75)--(-1,0) --cycle;
    \draw[axis] (-1,0)  -- (1,0) node[above, xshift=-1.5cm, yshift = 1cm] {$\bm{M_1(\alpha_1,e_1)}$} node[below, xshift=-1.5cm,yshift=-1cm] {$\bm{M_2(\alpha_2,e_2)}$};
  \end{tikzpicture}
  
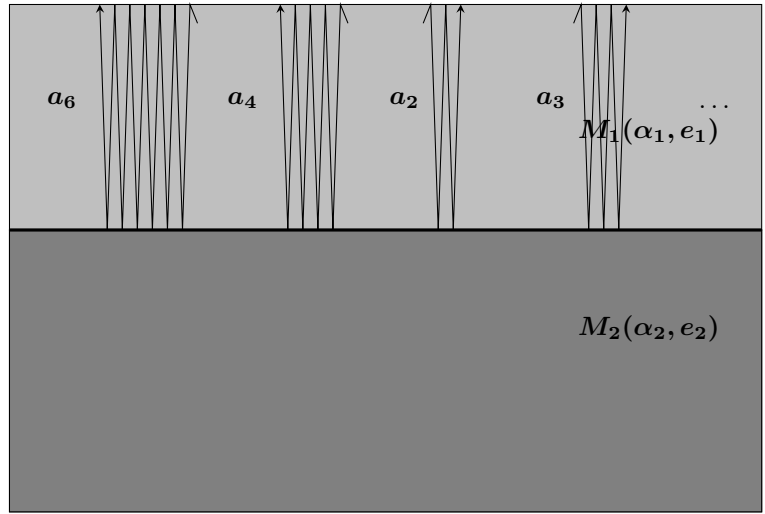
\captionof{figure}{Waves that are reflected first at the interface between $M_0$ and $M_1$}
  \end{center}
    Denote by $a_n$ the $n$\textsuperscript{th} reflection order wave, by $R_0(R_1)$ the reflection coefficient at the interface between $M_0$ and $M_1$($M_1$ and $M_2$, respectively). Let furthermore denote $T_0(T_0')$ the transmission coefficients at the interface between $M_0$ and $M_1$ in the upward (downward) direction. Since the single wave trains are reflected at the interface ``infinitely'' often, we obtain the following series corresponding to the surface temperature:
    \begin{align}\label{X}
    \begin{split}
        \sum_{n=0}^\infty a_n &= T_0A_0 \sum_{n=1}^\infty\left( 1 +  R_0^nR_1^n\exp{(-2n\sigma_1L_1)}\right)\\
        &=T_0A_0 \sum_{n=0}^\infty R_0^nR_1^n\exp{(-2n\sigma_1L_1)}\\
        &=T_0 A_0 \sum_{n=0}^\infty \left[R_0 R_1 \exp{(-2\sigma_1 L_1)}\right]^n\\
        &=T_0 A_0 \frac{1}{1-R_0R_1\exp{(-2\sigma_1L_1)}}\\
        &=\frac{T_0 A_0 }{1-R_0R_1\exp{(-2\sigma_1L_1)}},
        \end{split}
    \end{align}
    where we use the geometric series formula, since $|R_0 R_1 \exp{(-2\sigma_1 L_1)}| < 1$ holds true.
  
    \item[\textbf{(2)}] Waves that are first reflected from the interface between $M_1$ and $M_2$ as described by the following figure:
   \begin{center}
   \captionsetup{type=figure}
 \begin{tikzpicture}[
    scale=5,
    axis/.style={very thick, >=stealth'},
    important line/.style={thick},
    dashed line/.style={dashed, thin},
    pile/.style={thick, ->, >=stealth', shorten <=2pt, shorten
    >=2pt},
    every node/.style={color=black}
    ]
    \draw[fill=lightgray] (-1,0)--(-1,0.6)--(1,0.6)--(1,0)--(-1,0) --cycle;
    
    \draw[fill=gray] (-1,0)--(1,0)--(1,-0.75)--(-1,-0.75)--(-1,0) --cycle;
    
    \draw[->,>=stealth] (-0.5,0.55) -- (-0.52,0) -- (-0.54,0.6) -- (-0.56,0) -- (-0.58,0.6) -- (-0.6,0) -- (-0.62,0.6) node[above,xshift=-0.5cm,yshift=-1.5cm] {$\bm{b_3}$};
    \draw[->,>=stealth] (-0.1,0.55) -- (-0.12,0)  -- (-0.14,0.6) -- (-0.16,0) -- (-0.18,0.6) node[above,xshift=-0.5cm,yshift=-1.5cm] {$\bm{b_2}$};
    \draw[->,>=stealth] (0.5,0.55) -- (0.52,0) -- (0.54,0.6) -- (0.56,0) -- (0.58,0.6) -- (0.6,0) -- (0.62,0.6)--(0.64,0) -- (0.66,0.6) -- (0.68,0) -- (0.7,0.6) node[above,xshift=-1.3cm,yshift=-1.5cm] {$\bm{b_5}$} node[above, yshift=-1.5cm,xshift=1cm]{$\dots$};
    \draw[->,>=stealth] (0.1,0.55) -- (0.12,0) -- (0.14,0.6) -- (0.16,0) -- (0.18,0.6) -- (0.2,0) -- (0.22,0.6) node[above,xshift=-1cm,yshift=-1.5cm] {$\bm{b_3}$};
    \draw[axis] (-1,0)  -- (1,0) node[above, xshift=-1.5cm, yshift = 1cm] {$\bm{M_1(\alpha_1,e_1)}$} node[below, xshift=-1.5cm,yshift=-1cm] {$\bm{M_2(\alpha_2,e_2)}$};
  \end{tikzpicture}
  \captionof{figure}{Waves that are reflected first at the interface between $M_1$ and $M_2$}
  \end{center}    
    Denote by $b_n$ the $n$\textsuperscript{th} reflection order wave. Analogous to the first case, we obtain the following expression corresponding to the surface temperature fraction:
    \begin{align}
    \begin{split}
        \sum_{n=0}^\infty b_n &=T_0A_0R_1 \exp{(-2\sigma_1L_1)}\sum_{n=1}^\infty \left( 1 + R_0^nR_1^n\exp{(-2n\sigma_1 L_1)}\right)  \\
        &=T_0A_0R_1 \exp{(-2\sigma_1L_1)} \sum_{n=0}^\infty R_0^nR_1^n\exp{(-2n\sigma_1 L_1)}\\
        &= T_0A_0R_1 \exp{(-2\sigma_1L_1)} \sum_{n=0}^\infty \left[R_0R_1\exp{(-2\sigma_1L_1)}\right]^n\\
         &=T_0A_0R_1 \exp{(-2\sigma_1L_1)} \frac{1}{1-R_0R_1\exp{(-2\sigma_1L_1)}}\\
         &= \frac{T_0AR_1 \exp{(-2\sigma_1L_1)}}{1-R_0R_1\exp{(-2\sigma_1L_1)}},
         \end{split}
    \end{align}
    where we use the geometric series formula, since $|R_0R_1\exp{(-2\sigma_1L_1)}| < 1$ holds true.
\end{itemize}
By summing up both of series, we obtain the thermal wave interference expression at the surface:
\begin{align}\label{temperatur}
\begin{split}
\sum_{n=0}^\infty a_n + \sum_{n=0}^\infty b_n
&=T_0 A_0 \left[\frac{1+R_1 \exp{(-2\sigma_1 L_1)}}{1-R_0R_1\exp{\left(-2\sigma_1 L_1\right)}}\right]\eqqcolon \bm{\Tilde{T}}(x=0).
\end{split}
\end{align}
This provides us with the following expression for the time-dependent temperature at the surface: 
\begin{align}\label{temp}
    \bm{T}(x=0,t) = \Tilde{T}(x=0)\exp{\left[i(\omega t - \frac{\pi}{4})\right]}.
\end{align}
Notice that the phase shift in this formula, given by $\frac{\pi}{4}$, is excluded and has the meaning of normalizing the temperature by subtracting an ``infinitely'' thick layer,  cf. Equation \eqref{phaseshift}.

Since the wave vector $\sigma_1 = (1+i)\sqrt{\frac{\omega}{2\alpha_1}}=(1+i)\frac{1}{\mu_1}$ is complex, implying that the temperature amplitude is complex valued, we can split it into its real and its imaginary part. That provides us with an expression in polar coordinates. To this end we set 
\begin{align*}
    \sigma_1 = \text{Re}(\sigma_1) + i\text{Im}(\sigma_1) = \frac{1}{\mu_1} + i \frac{1}{\mu_1} \coloneqq \sigma_1' + i\sigma_1''.
\end{align*}
and calculate
\begin{align*}
\begin{split}
    &\text{Re}\left[\Tilde{T}(x=0)\right]\\
    &= \frac{1-R_0^2R_1^2\exp{(-4\frac{L_1}{\mu_1})}+R_1\exp{(-2\frac{L_1}{\mu_1})}\cos(-2\frac{L_1}{\mu_1})\left[1-R_1\right]}{\left[1-R_0^2R_1^2\exp{(-2\frac{L_1}{\mu_1})}\cos{(-2\frac{L_1}{\mu_1})}\right]^2 + \left[R_0R_1\exp{(-2\frac{L_1}{\mu_1})\sin{(-2\frac{L_1}{\mu_1})}}\right]^2}
    \end{split}
\end{align*}
and 
\begin{align*}
\begin{split}
    &\text{Im}\left[\Tilde{T}(x=0)\right]\\
    &= \frac{\left[ 1 + R_0\right]\left[R_1\exp{(-2\frac{L_1}{\mu})\sin{(-2\frac{L_1}{\mu_1})}}  \right]}{\left[1-R_0^2R_1^2\exp{(-2\frac{L_1}{\mu_1})}\cos{(-2\frac{L_1}{\mu_1})}\right]^2 + \left[R_0R_1\exp{(-2\frac{L_1}{\mu_1})\sin{(-2\frac{L_1}{\mu_1})}}\right]^2}.
    \end{split}
\end{align*}
Define $x\coloneqq \frac{-2L_1}{\mu_1}$. Then we obtain the amplitude $A_{\Tilde{T}}$ as
\begin{align}\label{amplitude}
\begin{split}
    \bm{A_{\Tilde{T}}} &=  \sqrt{\text{Re}\left[\Tilde{T}(x=0)\right]^2 + \text{Im}\left[\Tilde{T}(x=0)\right]^2}\\
    &=\frac{\sqrt{\left[1+R_1(1-R_0)\exp{(x)}\cos{(x)}-R_1^2R_0\exp(2x)  \right]^2+\left[R_1(1+R_0)\exp{(x)\sin{(x)}}  \right]^2}}{\left[1-R_0R_1\exp{(x)\cos{(x)}} \right]^2+\left[R_0R_1\exp{(x)\sin{(x)}} \right]^2}
    \end{split}
\end{align}
as well as the phase angle $\varphi_{\Tilde{T}}$
\begin{align}\label{phase}
\begin{split}
    \bm{\varphi_{\Tilde{T}}} &= \tan ^{-1}\left[\frac{\text{Im}\left[\Tilde{T}(x=0)\right]}{\text{Re}\left[\Tilde{T}(x=0)\right]}\right]\\
    &= \tan ^{-1}\left[\frac{\left[ 1 + R_0\right]\left[R_1\exp{(x)\sin{(x)}}  \right]}{1+ (1-R_0)R_1\exp{(x)}\cos{(x)}-R_1^2R_0\exp{(2x)}}\right]
    \end{split}
\end{align}
which provides us with an expression for the complex temperature amplitude at the surface in polar coordinates, i.e. 
$$\Tilde{T}(x=0)= A_{\Tilde{T}}\exp{(i\varphi_{\Tilde{T}})}.$$
The advantage of this form is that the amplitude as well as the phase are measurable real valued expressions. This is the reason, why these quantities (especially the phase angle) are used in the measuring process of layer thicknesses.

\subsubsection{Two-layered coatings on a substrate}

Let us now consider the following system of three layers $M_1$ - $M_3$ having thicknesses $L_1$ - $L_3$, where $M_3$ is assumed to be thermally thick: 

\begin{center}
\captionsetup{type=figure}
 \begin{tikzpicture}[
    scale=5,
    axis/.style={very thick, >=stealth'},
    important line/.style={thick},
    dashed line/.style={dashed, thin},
    pile/.style={thick, ->, >=stealth', shorten <=2pt, shorten
    >=2pt},
    every node/.style={color=black}
    ]
    \draw[fill=lightgray] (-1,0)--(-1,0.6)--(1,0.6)--(1,0)--(-1,0) --cycle;
     \draw[fill=LighterGray] (-1,0.6)--(-1,1.2)--(1,1.2)--(1,0.6)--(-1,0.6) --cycle;
    \draw[fill=gray] (-1,0)--(1,0)--(1,-0.75)--(-1,-0.75)--(-1,0) --cycle;
    
    \draw[->,>=stealth] (0,1.2) -- (0,1.7)node[above,yshift=1cm] {Transmitted thermal waves};
    \draw[->,>=stealth] (0.5,1.2) -- (0.5,1.7) ;
    \draw[->,>=stealth] (-0.5,1.2) -- (-0.5,1.7);
    \draw[->,>=stealth] (0,0) -- (0,-0.5)node[above,yshift=-1cm] {Transmitted thermal waves} ;
    \draw[->,>=stealth] (0.5,0) -- (0.5,-0.5) ;
    \draw[->,>=stealth] (-0.5,0) -- (-0.5,-0.5);

    \draw[very thick,<->,>=stealth] (1.1,0) -- (1.1,0.6) node[xshift=0.5cm,yshift=-1.5cm] {$\bm{L_{2}}$};
    \draw[very thick,<->,>=stealth] (1.1,0) -- (1.1,-0.75) node[xshift=0.5cm,yshift=2cm] {$\bm{L_3}$};
    \draw[very thick,<->,>=stealth] (1.1,0.6) -- (1.1,1.2) node[xshift=0.5cm,yshift=-1.5cm] {$\bm{L_{1}}$};

    \draw[cyan,->,>=stealth] (-0.5,1.2) -- (-0.52,0) -- (-0.54,0.6)  -- (-0.56,0) -- (-0.58,1.2) -- (-0.6,0) -- (-0.62,0.6) -- (-0.64,0) -- (-0.66,0.6) -- (-0.68,0) -- (-0.7,0.6) -- (-0.72,0) -- (-0.74,1.2);
    
    \draw[cyan,->,>=stealth] (-0.45,1.2) -- (-0.43,0.6) -- (-0.41,1.2);
    
    \draw[cyan,->,>=stealth] (-0.1,1.2) -- (-0.12,0.6)  -- (-0.14,1.2) -- (-0.16,0) -- (-0.18,1.2);

    \draw[cyan,->,>=stealth] (0.1,1.2)  -- (0.12,0.6) -- (0.14,0) -- (0.16,0.6) --(0.18,0) -- (0.2,1.2);
    
     \draw[cyan,->,>=stealth] (0.5,1.2) -- (0.52,0.6)  -- (0.54,0) -- (0.56,0.6) -- (0.58,0) -- (0.6,0.6) -- (0.62,0)  -- (0.64,1.2);
     \draw[axis] (-1,0)  -- (1,0) node[above, xshift=-1.5cm, yshift = 1cm] {$\bm{M_{2}(\alpha_2,e_2)}$} node[below, xshift=-1.5cm,yshift=-1cm] {$\bm{M_{3}(\alpha_3,e_3)}$} node[above, xshift=-1.5cm,yshift=4cm] {$\bm{M_{1}(\alpha_1,e_1)}$} node[above, xshift=-1.5cm,yshift=6.5cm] {$\bm{M_0(\alpha_0,e_0)}$};
  \end{tikzpicture}
  
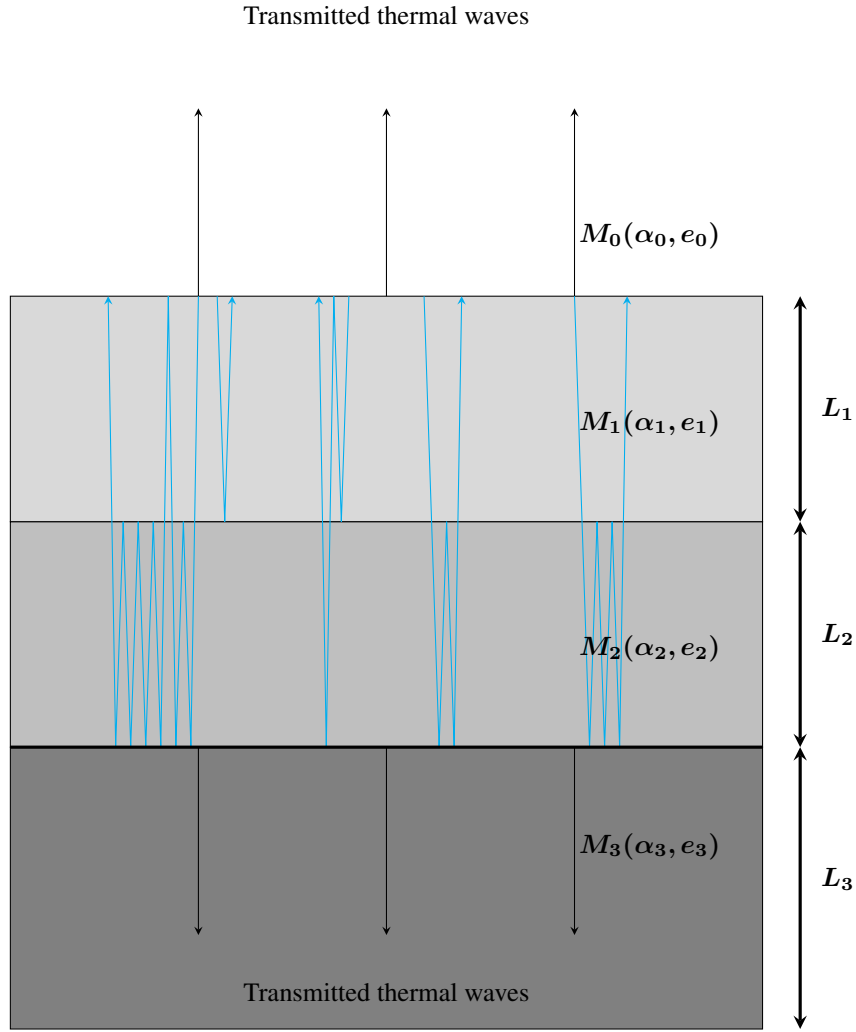
\captionof{figure}{System of two layers of coatings $M_1$ - $M_{2}$ and one layer of substrate $M_3$} \label{fig:ngleich2Fall}
  \end{center}

Obviously, adding a second coating layer increases the complexity in which thermal wave trains can possibly propagate in the system, see for example the blue arrows in Figure \ref{fig:ngleich2Fall}. The first interface for the case $n=1$ was a coating-to-substrate interface, while for $n=2$ the first interface is a coating-to-coating interface. Transmitted thermal waves can no longer be ignored, because the second layer is not thermally thick anymore, i.e. the reflections still contribute to the surface temperature significantly although they have a longer propagation path and are therefore damped stronger. A very elegant way of summarizing all possible wave trains is done by replacing a complex-valued \textit{effective reflection coefficient} $\Gamma_1$ for the real-valued reflection coefficient $R_1$ from Equation \eqref{temperatur}, cf. \cite{lepoutre1985nondestructive}, \cite{ALM1996}. 

Again, by introducing a new layer, the reflection process expands after passing the interface between $M_1$ and $M_2$ and takes also place at a lower level. Thus, interference effects occur in the new layer with material $M_2$ as well. This implies that the reflection coefficient $R_1$ does not suffice to describe the whole process anymore, so we need to change the expression in $R_1$ analogously to the derivation in Equation \eqref{X} in the following way:
\begin{align}
\begin{split}
    \bm{\Gamma_1} &\coloneqq  R_1 + T_1T_1'R_2\exp{(-2\sigma_2L_2)}\sum_{n=0}^\infty\left(R_1' R_2\exp{(-2\sigma_2L_2)}\right)^n\\
    &= R_1 + (1-R_1^2)R_2\exp{(-2\sigma_2L_2)}\sum_{n=0}^\infty\left(-R_1 R_2\exp{(-2\sigma_2L_2)}\right)^n\\
    &=R_1 + \left(R_2\exp{(-2\sigma_2L_2)}-R_1^2R_2\exp{(-2\sigma_2L_2)}\right)\frac{1}{1+R_1 R_2\exp{(-2\sigma_2L_2)}}\\
    &= \frac{R_1 + R_2\exp{(-2\sigma_2L_2)}}{1+R_1 R_2\exp{(-2\sigma_2L_2)}}.
    \end{split}
\end{align}
Here, $T_1$ ($T_1')$ denote the transmission coefficients in downward (upward) direction at the $M_1-M_2$ interface. By $R_1$ ($R_1'$) we denote the reflection coefficients at the same interface, where the thermal wave train stays in $M_1$ ($M_2$). Analogously, $R_2$ is the reflection coefficient for thermal waves in $M_2$ that stay in $M_2$.

Therefore, because of $R_1'=-R_1$, we use the geometric series formula again, what is possibl since $|-R_1R_2\exp{(-2\sigma_2L_2)}| < 1$.
Note that 
 \begin{align}
     \lim_{L_2\to \infty} \Gamma_1 = R_1,
 \end{align}
 i.e. we are effectively in the case $n=1$ again, if the second coating layer would be infinitely thick (or at least thermally thick).

We call the quantity $\Gamma_1$ the \textit{effective reflection coefficient}. Since the thermal interference processes are the same in the first layer as discussed before, we end up with the following surface temperature for the case $n=2$: 
\begin{align}\label{3}
\begin{split}
    \bm{T}(x=0,t) &= T_0 A_0 \left[\frac{1+\Gamma_1 \exp{(-2\sigma_1 L_1)}}{1-R_0\Gamma_1\exp{\left(-2\sigma_1 L_1\right)}}\right]\exp{\left[i\left(\omega t - \frac{\pi}{4}\right)\right]}\\
    &\eqqcolon \bm{\Tilde{T}}(x=0)\exp{\left[i\left(\omega t - \frac{\pi}{4}\right)\right]}.
    \end{split}
\end{align}

Once again, this expression provides us with a term for the amplitude,
\begin{align}
 \bm{A_{\Tilde{T}}} = \sqrt{\text{Re}\left[\Tilde{T}(x=0)\right]^2 + \text{Im}\left[\Tilde{T}(x=0)\right]^2},
    \end{align}
as well as a term for the phase angle,
\begin{align}
   \bm{\varphi_{\Tilde{T}}} = \tan ^{-1}\left[\frac{\text{Im}\left[\Tilde{T}(x=0)\right]}{\text{Re}\left[\Tilde{T}(x=0)\right]}\right],
\end{align}
yielding 
\begin{align*}
   \Tilde{T}(x=0)=A_{\Tilde{T}}\exp{(i\varphi_{\Tilde{T}})}.
\end{align*}
\newline 

\subsection{Generalization to multi-layered coating systems}

Subject of our investigations is to make use of the insights outlined before and generalize formula \eqref{3} for a multi-layered coating system, see Figure \ref{fig:mss}.
Consider the following system consisting of $n\in \mathbb{N}$ layers of coatings with material $M_1$ - $M_{n}$ having thicknesses $L_1$ - $L_n$ on a thermally thick substrate material $M_{n+1}$:
\begin{center}
\captionsetup{type=figure}
\begin{tikzpicture}[
    scale=5,baseline,
    axis/.style={very thick, >=stealth'},
    important line/.style={thick},
    dashed line/.style={dashed, thin},
    pile/.style={thick, ->, >=stealth', shorten <=2pt, shorten
    >=2pt},
    every node/.style={color=black}
    ]
    \draw[white] (-1,1.2)--(-1,1.4)--(1,1.4)--(1,1.2)--(-1,1.2) --cycle;
    \draw[fill=lightgray] (-1,0)--(-1,0.6)--(1,0.6)--(1,0)--(-1,0) --cycle;
     \draw[fill=LighterGray] (-1,0.6)--(-1,1.2)--(1,1.2)--(1,0.6)--(-1,0.6) --cycle;
      \draw[fill=LightererGray] (-1,1.4)--(-1,2)--(1,2)--(1,1.4)--(-1,1.4) --cycle;
    \draw[fill=gray] (-1,0)--(1,0)--(1,-0.75)--(-1,-0.75)--(-1,0) --cycle;
    
    \draw[->,>=stealth] (0,2) -- (0,2.4)node[above,yshift=1cm] {Transmitted thermal waves};
    \draw[->,>=stealth] (0.5,2) -- (0.5,2.4) ;
    \draw[->,>=stealth] (-0.5,2) -- (-0.5,2.4);
    \draw[->,>=stealth] (0,0) -- (0,-0.5)node[above,yshift=-1cm] {Transmitted thermal waves} ;
    \draw[->,>=stealth] (0.5,0) -- (0.5,-0.5) ;
    \draw[->,>=stealth] (-0.5,0) -- (-0.5,-0.5);
    \draw[very thick,<->,>=stealth] (1.1,0) -- (1.1,0.6) node[xshift=0.5cm,yshift=-1.5cm] {$\bm{L_{n}}$};
    \draw[very thick,<->,>=stealth] (1.1,0) -- (1.1,-0.75) node[xshift=0.5cm,yshift=2cm] {$\bm{L_{n+1}}$};
    \draw[very thick,<->,>=stealth] (1.1,0.6) -- (1.1,1.2) node[xshift=0.5cm,yshift=-1.5cm] {$\bm{L_{n-1}}$};
    \draw[very thick,<->,>=stealth] (1.1,1.4) -- (1.1,2) node[xshift=0.5cm,yshift=-1.5cm] {$\bm{L_{1}}$};

    \draw[cyan,->,>=stealth] (-0.5,1.2) -- (-0.52,0) -- (-0.54,0.6)  -- (-0.56,0) -- (-0.58,1.2) -- (-0.6,0) -- (-0.62,0.6) -- (-0.64,0) -- (-0.66,0.6) -- (-0.68,0) -- (-0.7,0.6) -- (-0.72,0) -- (-0.74,1.2);
    
    \draw[cyan, ->,>=stealth] (-0.45,1.2) -- (-0.43,0.6) -- (-0.41,1.2);
    
    \draw[cyan,->,>=stealth] (-0.1,1.2) -- (-0.12,0.6)  -- (-0.14,1.2) -- (-0.16,0) -- (-0.18,1.2);
    
   \draw[cyan,->,>=stealth] (-0.5,1.95) --(-0.5,1.4) ;
   
    \draw[cyan,->,>=stealth] (0.1,1.2)  -- (0.12,0.6) -- (0.14,0) -- (0.16,0.6) --(0.18,0) -- (0.2,1.2);
    
     \draw[cyan,->,>=stealth] (0.5,1.2) -- (0.52,0.6)  -- (0.54,0) -- (0.56,0.6) -- (0.58,0) -- (0.6,0.6) -- (0.62,0)  -- (0.64,1.2);
     
     \draw[cyan,->,>=stealth] (0.5, 1.95) --  (0.5,1.4);
   
    \draw[->,>=stealth] (-0.64,1.95) -- (-0.66,1.4) -- (-0.68,2) -- (-0.7,1.4) -- (-0.72,2) -- (-0.74,1.4) -- (-0.76,2);
    \draw[->,>=stealth] (-0.1,1.95) --  (-0.12,1.4) -- (-0.14,2) -- (-0.16,1.4) -- (-0.18,2);
    \draw[->,>=stealth] (0.1,1.95)  -- (0.12,2) -- (0.14,1.4) -- (0.16,2) --(0.18,1.4) -- (0.2,2);
    \draw[->,>=stealth] (0.6,1.95) -- (0.62,2)  -- (0.64,1.4) -- (0.66,2) -- (0.68,1.4) -- (0.7,2);
    \draw[->,>=stealth] (-0.45,1.95) -- (-0.45,2);
    \draw[->,>=stealth] (0.45,1.95) -- (0.45,2);
    \draw[axis] (-1,0)  -- (1,0) node[above, xshift=-2.5cm, yshift = 1cm] {$\bm{M_{n}(\alpha_{n},e_n)}$} node[below, xshift=-2.5cm,yshift=-1cm] {$\bm{M_{n+1}(\alpha_{n+1},e_{n+1})}$} node[above, xshift=-2.5cm,yshift=4cm] {$\bm{M_{n-1}(\alpha_{n-1},e_{n-1})}$} node[above, xshift=-2.5cm,yshift=8cm] {$\bm{M_1(\alpha_1,e_1)}$} node[above, xshift=-5cm,yshift=6.25cm] {$\bm{\vdots}$};
  \end{tikzpicture}
  
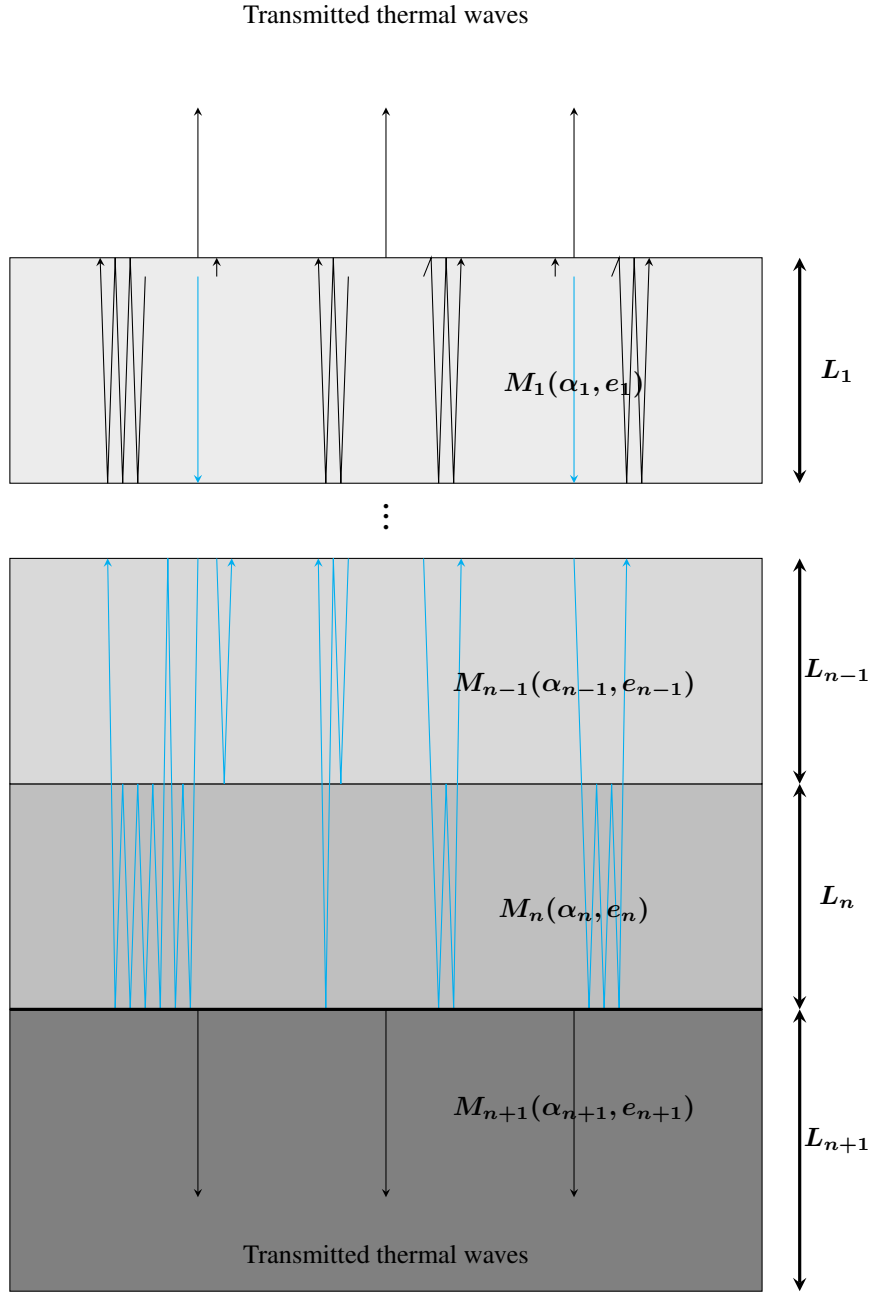
\captionof{figure}{System of $n$ layers of coating materials $M_1$ - $M_{n}$ and one layer of substrate material $M_{n+1}$} \label{fig:mss}
  \end{center}

For $n=2$, the idea was to look at the last reflection coefficient from the case $n=1$ and, since we add a new layer, replace it by some effective reflection coefficient $\Gamma_1$ including all potential wave trains that are apparent in the new system. This principle can be extended further, which leads to nested effective reflection coefficients, i.e. we define $\Gamma_j$ recursively by
\begin{align}
\begin{split}
&\Gamma_{n}\coloneqq R_{n},\\
    &\Gamma_{j}\coloneqq \frac{R_{j} + \Gamma_{j+1}\exp{(-2\sigma_{j+1}L_{j+1})}}{1+R_{j}\Gamma_{j+1}\exp{(-2\sigma_{j+1}L_{j+1})}},\ \ \  j=n-1,\dots,1.
    \end{split}
\end{align}

Although the generalized formula for the surface temperature of multi-layered coating systems looks as for the case $n=2$, i.e.

\begin{align}\label{finalFormula}
T(x=0,t) &= \underbrace{\left[ A_0 T_0\frac{1+\Gamma_1\exp{(-2\sigma_1L_1)}}{1-R_0\Gamma_1\exp{(-2\sigma_1L_1)}}\right]}_{=:A_{\Tilde{T}}\exp(i \varphi_{\Tilde{T}}),\ A_{\Tilde{T}}\in \mathbb{R}_+, \varphi_{\Tilde{T}} \in \mathbb{R}}\exp{\left[i\left(\omega t - \frac{\pi}{4} \right)\right]},
\end{align}

the effective reflection coefficient $\Gamma_1$ now contains all significant information about the coating system such as the layer thicknesses $(L_1,\dots,L_n)$, the thermal diffusivities $(\alpha_1,\dots,\alpha_n)$ and the thermal effusivities $(e_1,\dots,e_{n+1})$. Note that $\alpha_{n+1}$ and $L_{n+1}$ are not of any interest, because the substrate is thermally thick, but $e_{n+1}$ is needed for the calculation of $\Gamma_n$.

\section{Layer Thickness Determination}

Now that we have a mathematical model for the physical process of thermal wave interference, we want to formulate the inverse problem of determining the layer thicknesses $(L_1,\dots,L_n)$ from phase angle measurements. Also we want to address the issue of unknown thermal properties and present a solution in a separate subsection. 

\subsection{Forward and inverse operator}

As discussed in the preceding section, the thermal properties needed for the calculation of the surface temperature of a multi-layered coating system (with $n$ coating layers) can be summarized by the vector

\begin{align}
    \mathbf{p}:= (\alpha_1,\dots,\alpha_n,e_1,\dots,e_{n+1})^T \in \mathbb{R}_+^{2n+1}.
\end{align}

For $m\in \mathbb{N}$ and frequencies $\bm{\omega}:=(\omega_1,\dots, \omega_m)^T$, let us define the so-called \textit{forward operator}  

\begin{align}\label{forward_op}
 F_n\colon \mathbb{R}_+^{n}\to \mathbb{R}^m, \quad (L_1,\dots,L_{n})^T=:\bm{L}\mapsto \bm{\varphi_{\Tilde{T}}},
\end{align}

which maps the coating layer thicknesses $\mathbf{L}\in \mathbb{R}_+^n$ to the phase angle vector $\bm{\varphi_{\Tilde{T}}}\in\mathbb{R}^m$, i.e. component-wise for every frequency to the phase angle $\varphi_{\Tilde{T}}$ from Equation \eqref{finalFormula}. Note that in \eqref{forward_op} we dropped the dependencies on $\mathbf{p}$ and $\bm{\omega}$ since these are kept fixed, but it is sometimes practicable to use the extended notations,
\begin{align}
    \bm{\varphi}_{\Tilde{T}}= F_n(\mathbf{L}) = F_n(\mathbf{L}, \mathbf{p})=F_n(\mathbf{L}, \mathbf{p},\bm{\omega}).
\end{align}

For the determination of $(L_1,\dots,L_n)$, let us assume that an infrared camera measures the temperature response of such a system for every frequency. After some post-processing of the signals (cf. \cite{BUS1992}, \cite{wu1998lock},\cite{DIL2003}), we are given the phase angle data $$\bm{\varphi_{\Tilde{T},\textit{meas}}}\in \mathbb{R}^m.$$ 

We define the so-called \textit{inverse operator} by
    \begin{align*} \label{inverseOperator}
        G_n\colon \mathbb{R}^m\to \mathbb{R}_+^{n}
    \end{align*}
\begin{align}
    G_n(\bm{\varphi_{\Tilde{T},\textit{meas}}})\coloneqq \argmin\limits_{\bm{L}\in \mathbb{R}_+^{n}} ||F_n(\bm{L})-\bm{\varphi_{\Tilde{T},\textit{meas}}}||_2^2,
\end{align}

which maps the phase angle data to the layer thickness vector by minimizing a nonlinear least-squares functional. Here, $\|\cdot\|_2$ denotes the Euclidean norm.

\subsection{Unknown thermal properties and issues with all-at-once optimization}

For real-life applications the approach (Eq. \eqref{inverseOperator}) seems unusable as the thermal properties of every coating layer must be known, which usually is not the case. For example, in the automotive industry paint mixtures are changed occasionally and it does not make sense to perform a complex thermal analysis each time. In general, what is available for the calibration of any nondestructive testing device are different samples of a fixed multi-layered coating setup with varying but known coating layer thicknesses. These are typically measured either in a destructive (e.g. cross-sectional images under the microscope, etc.) or non-destructive (e.g. laser triangulation, x-ray, etc.) way. Whether this coating layer thickness information of every layer in every sample is gathered before or after the frequency scans and phase angle data collection with an infrared camera depends on the method in use. However, it is worth noting that this process only needs to be performed once for each coating system setup. 

In the following we present a concept to identify the needed thermal properties in order to determine $(L_1,\dots,L_n)$ of such coating systems.
As a first idea, we replace the functional from Equation \eqref{inverseOperator} by

\begin{align} \label{allatonceApproach}
    \min\limits_{(\bm{p},\bm{L})}||F_n(\bm{p},\bm{L})-\bm{\varphi_{\Tilde{T},\textit{meas}}}||_2^2,
\end{align}

which represents an all-at-once approach, i.e. all unknown parameters are determined simultaneously. Unfortunately, this approach shows the following ambiguity issue: For every $j=1,\dots,n$, terms of the form 
 \begin{align*}
     \exp(-2\sigma_j L_j) = \exp(-2(1+i)\sqrt{\frac{\omega}{2\alpha_j}}L_j)
 \end{align*}
imply
    \begin{align*}
        \frac{L_j}{\sqrt{\alpha_j}} = \frac{\Tilde{c} \cdot L_j}{\sqrt{\Tilde{c}^2 \cdot \alpha_j}} \text{ for } \Tilde{c} >0,
    \end{align*}

i.e. 
$$||F_n(\bm{p},\bm{L})-\bm{\varphi_{\Tilde{T},\textit{meas}}}||_2 = ||F_n(\Tilde{c}^2\bm{p},\Tilde{c}\bm{L})-\bm{\varphi_{\Tilde{T},\textit{meas}}}||_2,$$
since the thermal contrast between the layers and thus the reflection coefficients are identical. Hence, it is necessary to decouple the determination of $\mathbf{p}\in \mathbb{R}_+^{2n+1}$ and $\mathbf{L}\in \mathbb{R}_+^n$. This can be done by utilizing the sample data in a clever way.

Let us assume that $k\in \mathbb{N}$ samples ($n$ coating layers with different but known thicknesses $\bm{L}_j:=(L_{1,j},\dots,L_{n,j})^T$ for $j=1,\dots,k$) are available for the calibration process. Modulation of every sample surface by the frequencies $\bm{\omega}\in \mathbb{R}_+^m$ leads to a total of $m\times k $ phase angle values, i.e.

\begin{align*}
    \left(\bm{\varphi_{\Tilde{T},\textit{meas},\textit{j}}}\right)_{j=1,\dots,k}.
\end{align*}

For $1<k_1<k$, we divide the entire batch of sample data in two sets
\begin{align}
    \underbrace{S_1 = \left(\bm{\varphi_{\Tilde{T},\textit{meas},\textit{j}}}\right)_{j=1,\dots,k_1}}_{\text{Training or calibration set}} \text{ and } \underbrace{S_2 = \left(\bm{\varphi_{\Tilde{T},\textit{meas},\textit{j}}}\right)_{j=k_1+1,\dots,k}}_{\text{Test or confirmation set}}
\end{align}
and proceed with the following steps.

 STEP 1: Determine thermal properties with set $S_1$, i.e.
    \begin{align*}
    \Bar{\bm{p}} : = \argmin\limits_{\bm{p} \in \mathbb{R}_+^{2n+1}} \sum_{j=1}^{k_1}||F_n(\bm{p},\bm{L}_j)-\bm{\varphi_{\Tilde{T},\textit{meas},\textit{j}}}||_2^2.
\end{align*}

 STEP 2: Feed the thermal properties from STEP 1 into the objective functional and determine the coating layer thicknesses of sample $j$ by
    \begin{align*}
    \Bar{\bm{L}}_j : = \argmin\limits_{\bm{L} \in \mathbb{R}_+^{n}} ||F_n(\Bar{\bm{p}},\bm{L})-\bm{\varphi_{\Tilde{T},\textit{meas},\textit{j}}}||_2^2
\end{align*}
for $j=k_1+1,\dots,k$ . Evaluate the results  for set $S_2$ by calculating the error
\begin{align} \label{error}
    \sum_{j=k_1+1}^k \|\Bar{\bm{L}}_j-\bm{L}_j\|_2 ^2.
\end{align}

If the error \eqref{error} is sufficiently small, what of course depends on the accuracy requirements of the manufacturing process itself, the coating system setup is calibrated and tested successfully. If a higher accuracy is needed, further samples need to be processed or $k_1$ needs to be adjusted. Of course, the usage of more frequencies also improves the quality of the data.

\section{Conclusion}

In this article multi-layered coating systems have been investigated consisting of $n\in \mathbb{N}$ coating layers on a thermally thick substrate, which are periodically illuminated by a planar, sinusoidal wave form with a fixed frequency. This illumination generates a thermal wave with the same frequency, which is reflected and transmitted at layer interfaces. The surface temperature, which can be measured by an infrared camera, is a result of the superposition of all thermal wave trains propagating through the system. We developed a new model that describes the physical process of 1D thermal wave interference in such setups. This model describes the dependencies of the coating layer thicknesses, the frequency used and the thermal properties of the layers to the measured phase angle data. Given measured phase angle data, we then defined the inverse operator for computing the coating layer thicknesses. We also discussed the problem of unknown thermal properties and proposed a concept to determine these in advance. 

\section{Acknowledgments}

This research was funded by the European Fund for Regional Development from the Operational Program EFRE Saarland 2014-2020 with the objective "Investments in Growth and Employment".

\begin{minipage}{.5\textwidth}
  \centering
  \includegraphics[width=.8\linewidth]{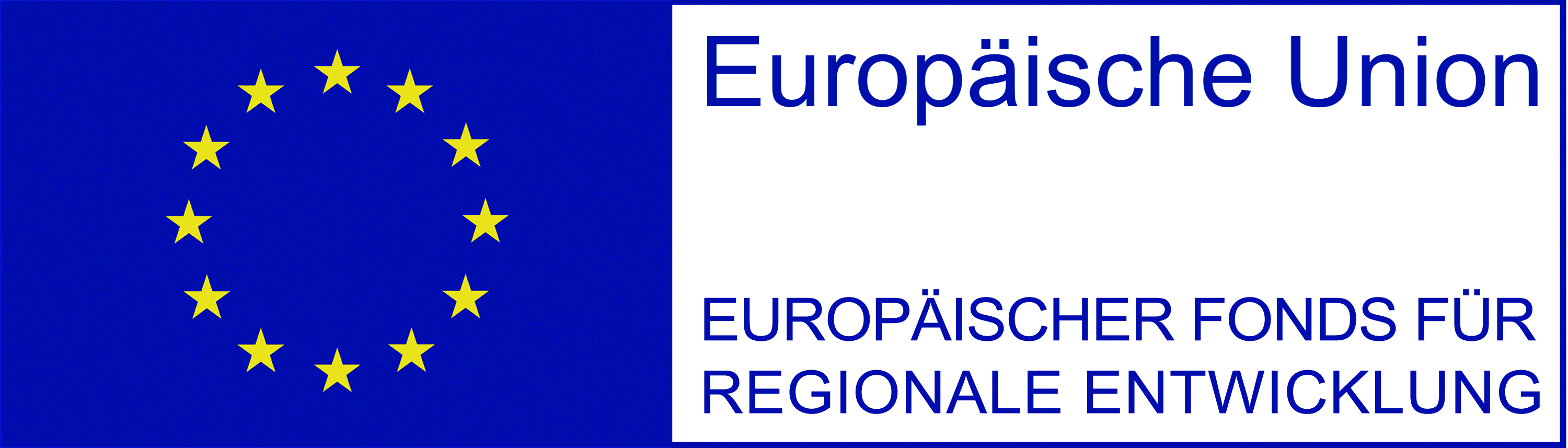}
\end{minipage}%
\begin{minipage}{.5\textwidth}
  \centering
  \includegraphics[width=.8\linewidth]{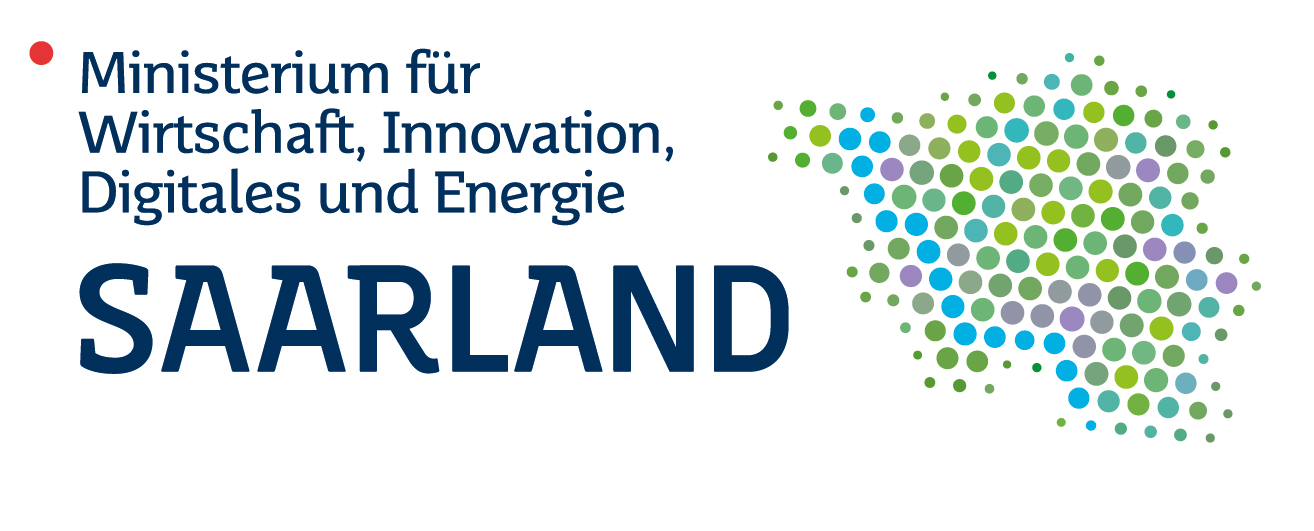}
\end{minipage}

\bibliography{references} 

\end{document}